\documentclass[12pt]{amsart}

\usepackage[utf8]{inputenc}

\usepackage{amsfonts}
\usepackage{amssymb}
\usepackage{amsmath}
\usepackage{amsthm}
\usepackage[T2A]{fontenc}

\usepackage{amsfonts}
\usepackage{amssymb}
\usepackage{amsmath}
\usepackage{amsthm}

\usepackage{amssymb}

\usepackage{enumitem}

\usepackage{graphicx}
\usepackage[T1]{fontenc}

\usepackage[utf8]{inputenc}

\usepackage[T2A]{fontenc}

\makeatletter
\@namedef{subjclassname@2010}{%
  \textup{2010} Mathematics Subject Classification}
\makeatother

\usepackage[T1]{fontenc}
\theoremstyle{definition}

\numberwithin{equation}{section}

\title[Best approximation vectors]{Best approximation vectors for completely irrational subspaces}

\author[N. Moshchevitin]{Nikolay Moshchevitin }
\address{
Israel Institute of Technology (Technion)}

\email{ moshchevitin@gmail.com, moshchevitin@technion.ac.il}

\date{}



\begin{document}
\maketitle

\begin{abstract}
  We give an easy  optimal bound for the dimension of the subspaces generated by the best Diophantine approximations.
\end{abstract}

 \section{Completely irrational subspaces.}

Let $m,n$ be positive integers and $d=m+n$. We consider Euclidean space $\mathbb{R}^d$ with  coordinates
$$
\pmb{z} = (\pmb{x},\pmb{y})^\top,\,\,\, \pmb{x} = (x_1,...,x_m)^\top \in \mathbb{R}^n,\,\,\, \pmb{y} = (y_1,...,y_n)^\top \in \mathbb{R}^n
$$
(for any martix $W$ by $W^\top$ we denote the matrix transposed to $W$).
Linear subspace $\mathcal{R} \subset \mathbb{R}^d$  of dimension $ r = {\rm dim}\, \mathcal{R}$ is called {\it rational} if the intersection
$\mathcal{R}\cap \mathbb{Z}^d$ is a $r$-dimensional sublattice in $\mathbb{Z}^d$. By $ H(\mathcal{R})$ we denote the {\it height} 
of the rational subspace $\mathcal{R}$, that is the fundamental volume of $r$-dimensional lattice $\mathcal{R}\cap \mathbb{Z}^d$.
We define  $m$-dimensional linear subspace $\mathcal{L}$ in $\mathbb{R}^d$ to be
{\it completely irrational} if for any rational subspace $\mathcal{R}\subset\mathbb{R}^d$  of dimension ${\rm dim }\, \mathcal{R}= n$
we have $\mathcal{L}\cap \mathcal{R} =\{ \pmb{0}\}$, that is $\mathcal{L}$ has no non-trivial intersections with rational subspaces of dimension $\le r$.

Let
$$\Theta
=\left(
\begin{array}{ccc}
\theta_{1,1}&\cdots&\theta_{1,m}\cr
\theta_{2,1}&\cdots&\theta_{2,m}\cr
\cdots &\cdots &\cdots \cr
\theta_{n,1}&\cdots&\theta_{n,m}
\end{array}
\right)\,\,\,
$$
be $m\times n$ real matrix. In the present paper we will consider only  $m$-dimensional subspaces of the form
\begin{equation}\label{ll}
\mathcal{L}_\Theta = \{ \pmb{z} \in \mathbb{R}^d: \,\,\,\, \pmb{y} = \Theta \pmb{x}\}.
\end{equation}
In the case of completely irrational subspace $\mathcal{L}_\Theta$ we will also call matrix $\Theta$ to be completely irrational.

\section{  Best approximation vectors.}

We use notation $|\pmb{\xi}| = \max_{1\le  j \le k} |\xi_j|$ for the sup-norm of the vector $\pmb{\xi}=(\xi_1,...,\xi_k)^\top\in \mathbb{R}^k$
and $||\pmb{\xi}|| = \min_{\pmb{a}\in \mathbb{Z}^k}|\pmb{\xi}-\pmb{a}|$ for the distance to the nearest integer point in sup-norm.

An integer vector $\pmb{x} $  is called a best approximation for the matrix $\Theta$ 
if
$$
||\Theta \pmb{x}|| =|\Theta\pmb{x} - \pmb{y}| <\min_{\pmb{x}'} ||\Theta \pmb{x}'||,
$$
where the minimum is taken over all non-zero integer points
$\pmb{x}'$ with $|\pmb{x}'|\le |\pmb{x}|$, $ \pmb{x}'\neq \pm\pmb{x}$,
and 
$ \pmb{y} = (y_1,...,y_n)^\top\in \mathbb{Z}^n$ is just the integer point at which the distance to the nearest integer point is attached.
For a best approximation $\pmb{x}$ satisfying this definition the point $-\pmb{x}$  will also be a best approximation.  We will  also consider the extended vector of the best approximation
$\pmb{z} = (\pmb{x},\pmb{y})^\top$.

We should 
note that in general it may happen that for two integer points $\pmb{x}'$ and  $\pmb{x}''$
with the same value 
$|\pmb{x}'|=|\pmb{x}''|
$
one has
$
 ||\Theta \pmb{x}'|| =  ||\Theta \pmb{x}''||$.
 So, in the general situation we cannot define uniquely a sequence of the best approximation vectors $\pm \pmb{x}_\nu\in \mathbb{Z}^m, \nu=1,2,3,...$  to satisfy
 \begin{equation}\label{1}
 |\pmb{x}_1|< |\pmb{x}_2|<...< |\pmb{x}_\nu|< |\pmb{x}_{\nu+1}|<....$$$$
 ||\Theta \pmb{x}_1|| >||\Theta \pmb{x}_2||>... > ||\Theta \pmb{x}_\nu||>... >||\Theta \pmb{x}_{\nu+1}||>... \, .   
 \end{equation}
 Also in the general situation we cannot define uniquely vectors $\pmb{y}_\nu\in \mathbb{Z}^n$
 to satisfy $ ||\Theta \pmb{x}_\nu|| =|\Theta\pmb{x}_\nu - \pmb{y}_\nu| $.
 However the  corresponding values of norms
 $ |\pmb{x}_\nu|$ and $ ||\Theta \pmb{x}_\nu||$ satisfying (\ref{1}) are well defined.
 We define matrix $\Theta$ to be {\it good}  if the unique sequence of
 best approximation vectors  $\pm\{\pmb{z}_\nu\}$  satisfying (\ref{1})  is well defined.   
 In this note we consider only good matrices.

\vskip+0.3cm

 Next, we should note that  parallelepiped
 \begin{equation}\label{pi}
 \Pi_\nu =\{ \pmb{z}= (\pmb{x},\pmb{y})^\top \in \mathbb{R}^d:\,\,\,\,\,
 |\pmb{x}|\le |\pmb{x}_{\nu+1}|,\,\,\,\,\,
 |\Theta\pmb{x} - \pmb{y}|\le ||\Theta\pmb{x}_\nu || \,\}
 \end{equation}
 has no non-zero integer points inside and 
 so  from Minkowski convex body theorem one has
 \begin{equation}\label{ppp}
 ||\Theta \pmb{x}_\nu|| \le  |\pmb{x}_{\nu+1}|^{-\frac{m}{n}}.
 \end{equation}
 Also we would like to recall the definition of irrationality measure function
 $$
 \psi_\Theta (t) = \min_{\pmb{x}\in \mathbb{Z}^m: \, 0< |\pmb{x}|\le t }\,\, ||\Theta\pmb{x}||,
 $$
 which is a piecewise constant function which is not continuous just in the points  $t_\nu = |\pmb{x}_\nu|$, as well as a general statement by 
 Jarn\'{\i}k \cite{J59} which claims the following.

 \vskip+0.3cm

(a) For any integers $   m \ge 2, n\ge 1$  and for any  function $ \psi (t)$ decreasing to zero as $t\to \infty$ there exist $n\times m$ matrices $\Theta$
with algebraically independent elements $\theta_{j,i} , 1\le i \le m, 1\le j \le n$ such that
\begin{equation}\label{z1}
 \psi_\Theta (t) \le \psi (t) \,\,\,\,\,\text{for all}\,\,\,t \,\,\,\text{large enough}.
\end{equation}

(b) For $ m =1$  and for any $n\ge 2$ the following holds.  Suppose that function $ \psi (t)$ decreas to zero as $t\to \infty$  but
$\lim_{t\to \infty} t\cdot \psi(t) = +\infty$. Then 
there exist $n\times 1$ matrices (vectors)  $\Theta$
with algebraically independent elements $\theta_{j,1},  1\le j \le n$ such that (\ref{z1}) holds.

\vskip+0.3cm

In particular,   Jarn\'{\i}k's result means that matrices $\Theta$ satisfying (\ref{z1})  can be chosen to be completely irrational and good.

 \section{  Dimension of subspaces of best approximations.}

For a good matrix $\Theta$ we define the value
$$
R(\Theta) = \min \{ s\in \mathbb{Z}_+:\,\,\, \text{there exists  a linear subspace }\,\,\mathcal {R} \subset \mathbb{R}^d\,\,
\text{of dimension}
$$
$$   \,\,\,\,\,\,\,\,\,\,\,\,\,\,\,\,\,\,\,\,\,\,\,\,\,\,\,\,\,\,\,\,\,\,\,\,\,\,\,\,\,\,\,{\rm dim}\, \mathcal{R} =s \,\,\,
\text{and}\,\, \nu_0 \in \mathbb{Z}_+
\,\,\,
\text{such that }\,\,\,
\pmb{z}_\nu \in \mathcal{R}\,\,\,\text{for all}\,\,\, \nu\ge \nu_0\}
$$
which will be the main object of our study.   It is clear that for any $m,n\ge 1$ for any good matrix 
$\Theta$ 
one has $ R(\Theta) \ge 2$.
In the next theorem we collect together all the results 
  about all possible values of $R(\Theta)$ for completely irrational $\Theta$.

\vskip+0.3cm

{\bf Theorem 1.}
{\it Suppose $\Theta$ to be a completely irrational good matrix. Then 

\vskip+0.2cm

{\rm (1)} if $ m=1$ for any $n$ one has $ R(\Theta) = d= n+1$;

{\rm (2a)} if $n=1$ and $ m \ge 2$ one has  $ R(\Theta) \ge 3$;

{\rm (2b)} if $n=1$  for any  $ m \ge 2$ there exist completely irrational good $\Theta$ with    $ R(\Theta) = 3$;

{\rm (3a)} if $n,m\ge 2 $  and  $ n\le m $ one has    $ R(\Theta) \ge n+2$;

{\rm (3b)} if $n,m\ge 2 $  and  $ n\le m $ there exist completely irrational  good $\Theta$ with    $ R(\Theta) = n+2$;

{\rm (4a)} if $n,m\ge 2 $  and  $ n> m $ one has    $ R(\Theta) \ge n+1$;

{\rm (4b)} if $n,m\ge 2 $  and  $ n> m $ there exist completely irrational  good $\Theta$ with    $ R(\Theta) = n+1$.

}

\vskip+0.3cm

Statements (1), (2a), (2b) are well known and discussed in surveys \cite{Che,Msing}.
We should note that 
statement (2a) goes back to Jarn\'{\i}k \cite{Jsze,Jche} and was used by Davenport and Schmidt \cite{DS}, Lemma 5.
Various examples of degeneracy of dimension from \cite{Msing} deal with the situation when $\Theta$ is not a completely irrational matrix.
Statement (2b) was first proven in  \cite{MDAN}. In \cite{Neck} (Theorem 3.5) it is shown that 
for an $m$-dimensional completely irrational subspace the inequality $ R(\Theta) \ge n+1$ is valid (see also Corollary 2 of Theorem 9 from \cite{Msing}), 
and so (4a) follows from this fact.
 Statements (3a), (3b), (4b) has never been documented. We give a proof of statement (3a) in Section 4 below. Statement (4b) is proven in Section 5. The proof of statement (3b) is quite similar and we give a sketched proof in Section 6.
 
 \vskip+0.3cm
 We would like to formulate here some additional comments concerning cases  (3) and (4). To do this we need to recall the definition of {\it uniform Diophantine exponent}
 $$\hat{\omega} (\Theta) = 
 \sup\{ \gamma \in \mathbb{R}:\,\,\, \limsup_{t\to \infty }\,\, t^\gamma \cdot \psi_\Theta ( t) =\limsup_{\nu \to \infty}\,\, |\pmb{x}_{\nu+1}|^\gamma \cdot  ||\Theta\pmb{x}_\nu ||<+ \infty\}.
 $$
 It is well known that $\hat{\omega} (\Theta)$ satisfies
 $$
 \frac{m}{n} \le \hat{\omega}(\Theta) \le 
 \begin{cases} 1,\,\,\,\,\, \,\,\,\,\, m=1\cr
 \infty,\,\,\,\,\, \,\,m \ge 2
 \end{cases}
 $$ 
 and may attain any value in these intervals.
 
 \vskip+1cm

{\bf Remark 2.}

(a) Suppose that $  n>  m\ge 2$ and $ R(\Theta) = n+1$. Then $ \frac{m}{n} \le \hat{\omega} (\Theta)\le 1$. Moreover, $\hat{\omega}  (\Theta)$ may attain
 any value in this interval;

 (b) Suppose that $  2\le n\le m$ and $ R(\Theta) = n+2$. Then, of course,  inequality $ \frac{m}{n} \le \hat{\omega}  (\Theta)\le +\infty $ holds, and moreover, $\hat{\omega}  (\Theta)$ may attain
 any value in this interval.
 
 (c) In particular, there exist matrices $\Theta$ with degenerate dimension of subspaces of best approximation vectors and with  $\hat{\omega} (\Theta)=\frac{m}{n}$.

 \vskip+0.3cm
Remark 2 follows form the fact that irrationality measure functions $\psi_\Theta (t)$ and $\psi_\Theta (t)$ eventually coincide (see Lemma 5 below) and the irrationality measure function defines the uniform exponent.

\vskip+0.3cm
Quite recently Schleischitz  \cite{Sch} studied various properties of the sets of $\Theta$ for which
$R(\Theta) <d =m+1$ in the case $n=1$. In particular he give  various bounds for Hausdorff dimension of the sets under the consideration. The analysis of the subspaces of best approximations and their dimensions turned out to be important to some other studies related to  Diophantine exponents, see \cite{Sch0}. In particular, in  \cite{Sch0} it is shown that  for the
{\it ordinary Diophantine exponent}
 $${\omega} (\Theta) \! = \!
 \sup\{ \gamma \in \mathbb{R}:\liminf_{t\to \infty } t^\gamma \cdot \psi_\Theta ( t) =\limsup_{\nu \to \infty}\, |\pmb{x}_{\nu+1}|^\gamma \cdot  ||\Theta\pmb{x}_\nu ||<+ \infty\}\!
 \ge \!\hat{\omega}(\Theta)
 $$
 for $m\times n $ matrix $\Theta$ 
 under the condition  $  r= R(\Theta) <d$ one has non-trivial lower bounds of the form
 $$
 {\omega} (\Theta) \ge G(m,n,r) \cdot \frac{m}{n},\,\,\,\,\text{where}\,\,\,\,\,
 G(m,n,r) >1.
 $$
 This happens due to existence of a non-trivial lower bound for the ratio $\frac{\omega(\Theta)}{\hat{\omega}(\Theta)}$.

\vskip+0.3cm
If we suppose that $\Theta$ is not necessary completely irrational but good, certain existence results were obtained in \cite{Msing}, Section 2. In particular (Theorems 7, 10, 13, 14), it was shown that 

(1) for any $m\ge 2$ and $ n >m$ there exist good matrices $\Theta$ with entries $\theta_{i,j}$ linearly independent together with 1 over $\mathbb{Q}$ such that  $ R(\Theta) = 2$;

(2) for any $ m>n$ for a good matrix $\Theta$ one has $ R(\Theta) \ge 3$;

(3) for any $m\ge 2$ and arbitrary $ n\ge 1$ there exist good matrices $\Theta$ with entries $\theta_{i,j}$ linearly independent together with 1 over $\mathbb{Q}$ such that  $ R(\Theta) = 3$.

\vskip+0.3cm
We finish this section wth the following remark. In the case $m=n=2, d=4$ we do not know if there exists good matrix $\Theta$ with $R(\Theta)=2$. By our Theorem 1 (3a) we see that this  is not possible for a completely irrational  matrix $\Theta$, as in such a case one has $R(\Theta) = 4$
(see also Corollary 3.7 from \cite{Neck}). In  \cite{Msing} it is shown that  $R(\Theta)=3$ is not possible for good  $2\times 2$ matrices.

\section{ Proof of the statement (3a).}

We assume that $\mathcal{L}$ defined in (\ref{ll}) is a completely irrational subspace. 
As $R(\Theta)\ge n+1$  (see Theorem 3.5 from \cite{Neck})
it is enough to prove  that $R(\Theta)= n+1$   does not hold.
Suppose that $\mathcal{R} \subset \mathbb{R}^d, {\rm dim}\, \mathcal{R} = n+1$ be the rational subspace of the smallest dimension form the definition of $R(\Theta)$.

We consider the intersection 
$ \ell = \mathcal{R}\cap \mathcal{L}$. We claim that $\ell$ is a one-dimensional linear subspace. Indeed,
if ${\dim}\, \ell \ge 2$ then the union
$\mathcal{R}\cup \mathcal{L} $ is contained in a $(d-1)$-dimensional subspace $\mathcal{W}$ of $\mathbb{R}^d$. Then for any $n$-dimensional   (rational) subspace $ \mathcal{R}'\subset \mathcal{R}$,
both $n$-dimensional subspace $\mathcal{R}'$ and $m$-dimensional subspace $\mathcal{L}$ lie 
in $(m+n-1)$-dimensional subspace $\mathcal{W}$. So for rational subspace $\mathcal{R}'$ we have
$\mathcal{R}'\cap \mathcal{L} \neq \{ \pmb{0}\}$ and this is not possible by the definition of completely irrational subspace.

As $ n\le m$ from (\ref{ppp}) we see that 
\begin{equation}\label{pp}
 ||\Theta \pmb{x}_\nu||  \cdot  |\pmb{x}_{\nu+1}|\le 1.
\end{equation}
We deal with two successive best approximations $ \pmb{z}_j = (\pmb{x}_j,\pmb{y}_j)^\top\in \mathcal{R}, j = \nu, \nu+1$. We prove that 
\begin{equation}\label{pp}
\lim_{\nu \to \infty }
 ||\Theta \pmb{x}_\nu||  \cdot  |\pmb{x}_{\nu+1}| =\infty
\end{equation}
and this will contradict  to (\ref{pp}).

Let  
${\rm dist}\, (\mathcal{A},\mathcal{B}) $ denotes the Euclidean distance between the sets
 $\mathcal{A},\mathcal{B}\subset \mathbb{R}^d$ and 
${\rm angle}\, (\pmb{z}',\pmb{z}'')$ denote the angle between non-zero vectors in $\mathbb{R}^d$.
We consider the orthogonal complement $\ell^\perp$ which is a $(d-1)$-dimensional subspace and subspaces
$\frak{R} = \mathcal{R}\cap \ell^\perp, \frak{L} = \mathcal{L}\cap \ell^\perp$.
By compactness argument we see that 
$$
\delta (\mathcal{R},\mathcal{L} ) =
\min_{\pmb{z}' \in \frak{R},\pmb{z}''\in \frak{L}} \, {\rm angle}\, (\pmb{z}',\pmb{z}'') >0.
$$
It is clear that for any $\pmb{z}  = (\pmb{x},\pmb{y})^\top \in \mathcal{R}$ one has
\begin{equation}\label{qq}
{\rm dist}\, (\pmb{z}, \ell) \le \frac{{\rm dist}\, (\pmb{z},\mathcal{L})}{\sin \delta (\mathcal{R},\mathcal{L} )}
\le C_{\mathcal{L},\mathcal{R} }\, |\Theta\pmb{x}-\pmb{y}|
\end{equation}
with some positive  constant $C_{\mathcal{L},\mathcal{R} }$
and
\begin{equation}\label{qqq}
{\rm dist}\, (\pmb{z}_\nu, \ell) 
\to 0,\,\,\,\,\,\nu \to \infty.
\end{equation}

Now we consider  two-dimensional subspace $\pi_\nu =\langle \pmb{z}_\nu, \pmb{z}_{\nu+1}\rangle_{\mathbb{R}}$ and two-dimensional lattice $\Lambda_\nu = \langle \pmb{z}_\nu, \pmb{z}_{\nu+1}\rangle_{\mathbb{Z}}\subset \pi_\nu $. Let $\Delta_\nu $ denotes covolume  (the area of the fundamental parallelogramm) of $\Lambda_\nu$. As the number of two-dimensional integer sublattices of $\mathbb{Z}^d$ with bounded covolume is finite and $\ell \cap \pi_\nu = \{\pmb{0}\} \,\, \forall\, \nu$, from (\ref{qqq}) we see that 
\begin{equation}\label{qqqq}
\Delta_\nu \to \infty, \nu \to \infty.
\end{equation}
As $ \pmb{z}_\nu, \pmb{z}_{\nu+1} \in \Pi_\nu \cap\mathbb{R} , \ell \subset \Pi_\nu \cap \mathbb{R}$ from
(\ref{qq}) we see that 
$$
\pm \pmb{z}_\nu, \pm\pmb{z}_{\nu+1} \in \pi_\nu\cap \Omega,
$$
{where}
$$
\Omega = 
\{ \pmb{z}\in \mathcal{R}: \,\,|\pmb{z}|\le C_1(\Theta)|\pmb{x}_{\nu+1}|, \,\,\,{\rm dist}\, (\pmb{z}, \ell) \le C_2 (\Theta, \mathcal{R}) ||\Theta\pmb{x}_\nu||
\}
$$
with some positive $ C_1(\Theta), C_2 (\Theta, \mathcal{R})$ and so
\begin{equation}\label{co}
\Delta_\nu \le \,\text{area of}\, \,\, \pi_\nu\cap \Omega \le   C_3 (\Theta, \mathcal{R})\, |\pmb{x}_{\nu+1}|\cdot
||\Theta\pmb{x}_\nu||,  \,\,\, C_3 (\Theta, \mathcal{R}).
\end{equation}
From (\ref{qqqq}) and (\ref{co}) we get (\ref{pp}).$\Box$.

\section{  Proof of the statement (4b).}

We consider vector  $\Theta^{*} = (\theta^*_1,...,\theta_n^*)^\top \in \mathbb{R}^n$ which consist of algebraically independent components $\theta^*_j$.
Let 
$$
\psi_{\Theta^*} (t) =   \min_{x\in \mathbb{Z}_+:\, x\le t}  ||\Theta^* x|| =
\min_{x\in \mathbb{Z}_+:\, x\le t}  \max_{1\le j \le n} ||\theta_j^* x||
$$
be the irrationality measure function for simultaneous approximation with vector $\Theta^1$
and $( x_\nu,y_{1,\nu},...,y_{n,\nu}) \in \mathbb{Z}^{n+1}, \nu =1,2,3,...$ be the sequence  of all best approximations to $\Theta^*$. It is clear that $ R(\Theta^*) = n+1$. Denote
$$
\xi_\nu =
\psi_{\Theta^*} (x_\nu) = ||\Theta^*x_\nu||.
$$
We consider $m\times n$ matrices $\Theta$   for the form
\begin{equation}\label{mma}
\Theta =(\Theta^*|\Theta^2),\,
\text{where}\,
\Theta^2 
=\left(
\begin{array}{ccc}
\theta_{1,2}&\cdots&\theta_{1,m}\cr
\theta_{2,2}&\cdots&\theta_{2,m}\cr
\cdots &\cdots &\cdots \cr
\theta_{n,2}&\cdots&\theta_{n,m}
\end{array}
\right)\,
\text{is a }\, (m-1)\times n\,
\text{matrix}.
\end{equation}
We identify the set of all matrices $\Theta^2$  with $\mathbb{R}^N, N = (m-1)n$ and consider Lebesgue measure on $\mathbb{R}^N$.
We consider the sequence of integer vectors
\begin{equation}\label{ss}
(x_\nu,\underbrace{0,...,0}_{(m-1)\,\,\text{times}}; y_{1,\nu},...,y_{n,\nu}), \,\,\, \nu =1,2,3,... \, .
\end{equation}

\vskip+0.3cm
{\bf Lemma 3.}\,{\it Assume that the series
\begin{equation}\label{ser}
\sum_{\nu=1}^\infty x_{\nu+1}^{m} \xi_\nu^n
\end{equation}
converges.
Then for almost all matrices $\Theta^2$  the sequence of the extended best approximation vectors for matrix $\Theta$  defined in (\ref{mma}) differs from the sequence 
(\ref{ss}) by at most finite number of elements. In particular, for almost all  $\Theta^2$   there exists $t_0$ such that
\begin{equation}\label{eeee}
\psi_\Theta (t) = \psi_{\Theta^*} (t),\,\,\,\,\,\, \forall\, t \ge t_0.
\end{equation} }
\vskip+0.3cm
Statement (4b) of Theorem 3 follows immediately from Lemma 5. Indeed, we see that  $R(\Theta) = R(\Theta^*) = n+1$.  Assuming $ n>m$, take $ \gamma \in  \left(\frac{m}{n},1\right)$. Then 
by Jarn\'{\i}k's result cited in Section 2, part (b) applied with $ \psi(t) = t^{-\gamma}$, we get $\Theta^*$ with algebraically independent elements and  $ \psi_{\Theta^*} (t)  \le t^{-\gamma}$. This means that for $\nu$ large enough one has
\begin{equation}\label{ppr}
x_{\nu+1}^m\xi_{\nu}^{n} \le x_{\nu+1}^{\gamma_1},\,\,\,\, \gamma_1 = m-n\gamma <0,
\end{equation}
and the series (\ref{ser}) converges as $x_\nu$ grow exponentially (see Lemma 1 from \cite{BL}). Now Lemma  5 ensures the existence of completely irrational extended matrix
$\Theta =(\Theta^*|\Theta^2)$ satisfying the desired properties.

\vskip+0.3cm
Proof of Lemma 3.
For vector $\pmb{x} = (x_1,...,x_m)\in \mathbb{Z}^m$ we define a shortened vector
$\underline{\pmb{x}}= (x_2,...,x_m)$.

To show that sequence (\ref{ss}) eventually coincides with the sequence of the best approximations for matrix $\Theta$ with entries $\theta_{i,j} \in [0,1]$ it is enough to show that there exists $\nu_0$ such that for all $\nu \ge \nu_0$ and for all
integer vectors
\begin{equation}\label{222}
(\pmb{x}, \pmb{y})
 = (x_1,x_2,...,x_m; y_1,...,y_n),
\in \mathbb{Z}^{d} 
$$$$
\text{with}
\,\,\,\, 0\neq
 |\underline{\pmb{x}}|\le |\pmb{x}|\le x_{\nu+1}
,
\,\,\,
\max_{1\le j\le n} |y_j - x_1\theta_j^*|
\le m  |\underline{\pmb{x}}|
\end{equation}
one has
$$
|\Theta \pmb{x} - \pmb{y}| > \xi_\nu
$$
(we should note that for fixed $\pmb{x}$ the number of integer vectors 
$\pmb{y}\in \mathbb{Z}^n$ satisfying the last inequality  in (\ref{222}) is $\ll_m |\underline{\pmb{x}}|^n$).
This condition may be rewritten as follows.    Define the sets
$$
\Omega_{j,\nu} ( \pmb{x} , y_j) =
\{ (
\theta_{j,2},...,\theta_{j,m}) \in [0,1]^{m-1}:\,
|\theta_{j,2}x_2+...+\theta_{j,m}x_m  + (\theta_{1}^* x_1-y_j)|\le \xi_\nu\}.
$$
As  in our consideration $|\underline{\pmb{x}}|\neq 0$, these sets 
$
\Omega_{j,\nu} ( \pmb{x} , y_j) $ are intersections of strips of thickness $\frac{\xi_\nu}{|\underline{\pmb{x}}|}$ around hyperplanes
$$
\{ (
w_{j,2},...,w_{j,m}) \in \mathbb{R}^{m-1}:\,
w_2x_2+...+w_mx_m  + (\theta_{1}^* x_1-y_j)= 0\}
$$
with the unit cube $[0,1]^{m-1}$. 
Consider
$$
\Omega_\nu( \pmb{x} , \pmb{y})  =
\Omega_{1,\nu} ( \pmb{x} , y_1)\times \cdots \times \Omega_{n,\nu}( \pmb{x} , y_n)   \subset [0,1]^N
$$
and define
$$
\frak{W}_\nu = 
\bigcup_{\pmb{x}
}
\bigcup_{\pmb{y}}
\Omega_\nu( \pmb{x} , \pmb{y})  
$$
where the union is taken over all vectors $\pmb{x},\pmb{y}$ satisfying (\ref{222}).
Now the condition
\begin{equation}\label{main}
\exists\, \nu_0:\,\,\,
\forall \, \nu \ge \nu_0 
\,\,\,\,\,
\text{one has}
\,\,\,\,\,
\Theta^2 \not \in \frak{W}_\nu
\end{equation}
ensures that  the 
sequence (\ref{ss}) eventually coincides with the sequence of the best approximations for matrix $\Theta$.

For the $(m-1)$-dimensional measure of $
\Omega_j ( \pmb{x} , y_j)$ we have an upper bound
\begin{equation}\label{bou}
\mu_{m-1} \left(
\Omega_j ( \pmb{x} , y_j)\right) \ll_m \frac{\xi_\nu}{|\underline{\pmb{x}}|},
\end{equation}
and so for the measure of $\frak{W}_\nu$ we get
$$
\mu_N (\frak{W}_\nu) \ll_{m} \sum_{\pmb{x}}
\sum_{\pmb{y}} \left(\frac{\xi_\nu}{|\underline{\pmb{x}}|}\right)^n\ll_{m,n} \xi_\nu^n \sum_{\pmb{x}} 1\ll_{m,n} \xi_\nu^n\cdot{x}_{\nu+1}^m.
$$
Now everything follows from Borel-Cantelli lemma.$\Box$

\section{ Sketched proof of the statement (3b).}

The proof of statement (3b) is similar to those of (4b). We give a sketch of a proof and left the details to the reader.
We may assume $ m\ge 3$.
One should consider matrices
\begin{equation}\label{mma1}
\Theta =(\Theta^*|\Theta^2),\,\,\,\,
\text{where}\,\,\,\,
\Theta^*
=\left(
\begin{array}{cc}
\theta_{1,1}^*&\theta_{1,2}^*\cr
\theta_{2,1}^*&\theta_{2,2}^*\cr
\cdots &\cdots \cr
\theta_{n,1}^*&\theta_{n,2}^*
\end{array}
\right),
\,\,\,\,
\Theta^2 
=\left(
\begin{array}{ccc}
\theta_{1,3}&\cdots&\theta_{1,m}\cr
\theta_{2,3}&\cdots&\theta_{2,m}\cr
\cdots &\cdots &\cdots \cr
\theta_{n,3}&\cdots&\theta_{n,m}
\end{array}
\right).
\end{equation}
Then we suppose $\Theta^*$ to be completely irrational.
By Theorem 1,  statement (3a) we see that $R(\Theta^*) = n+2$. Note that now
$$
\psi_{\Theta^*} (t) =   \min_{\pmb{x}\in \mathbb{Z}^2:0< |\pmb{x}|\le t}  ||\Theta^* \pmb{x}|| =
\min_{x_1,x_2\in \mathbb{Z}:\,0<\max(|x_1|,|x_2|) \le t}  \,\,\,  \max_{1\le j \le n}\,\, ||\theta_{j,1}^* x_1+\theta_{j,2}^*x_2||,
$$
and for the best approximation vectors $ \pmb{z}_\nu = (x_{1,\nu},x_{2,\nu};  y_{1,\nu},...,y_{n,\nu})\in \mathbb{Z}^{n+2}$ for $\Theta^*$ we consider
\begin{equation}\label{dde}
\xi_\nu =
\psi_{\Theta^*} (\pmb{x}_\nu),\,\,\, \pmb{x}_\nu = (x_{1,\nu},x_{2\nu}).
\end{equation}
Instead of Lemma 3 we need the following 
\vskip+0.3cm

{\bf Lemma 4.}\,{\it Assume that the series (\ref{ser}) for $\xi_\nu$ defined in (\ref{dde}) 
converges.
Then for almost all matrices $\Theta^2$  the sequence of the extended best approximation vectors for matrix $\Theta$  defined in (\ref{mma1}) differs from the sequence 
$$
(x_{1,\nu},x_{2,\nu},\underbrace{0,...,0}_{(m-2)\,\,\text{times}}; y_{1,\nu},...,y_{n,\nu}) \in \mathbb{Z}^{n+m}, \,\,\, \nu =1,2,3,... 
$$
associated to the sequence
$\pmb{z}_\nu \in \mathbb{Z}^{n+2}$ of
the best approximations for matrix $\Theta^*$ by at most finite number of elements, and in particular (\ref{eeee}) holds.
}
\vskip+0.3cm
The proof is quite similar to the proof of Lemma 5. It also uses Borel-Cantelli argument.

\vskip+0.3cm

To finish the proof of
statement (3b) of Theorem 1
one needs  
to apply 
 Jarn\'{\i}k's result cited in Section 2, part (b) applied with $ \psi(t) = t^{-\gamma}$,
 $\frac{m}{n} < \gamma < +\infty$.
 So one
gets  $\Theta^*$ with algebraically independent elements and  $ \psi_{\Theta^*} (t)  \le t^{-\gamma}$. Now again for  $\nu$ large enough one has
(\ref{ppr})
  the convergence of the series (\ref{ser}) follows from the exponential growth o $|\pmb{x}_\nu|$.

 \vskip+0.3cm
 {\bf Acknowledgements}.
 This work has received funding from the European Research Council (ERC) under the European Union’s Horizon 2020 Research and Innovation Program, Grant agreement no. 754475

 \vskip+1cm

 \end{document}